\documentclass[12pt]{article}
\usepackage{amssymb,amsmath,amsthm}
\numberwithin{equation}{section}

\makeatletter
\renewcommand{\@makefnmark}{\hbox{\mathsurround=0pt $^\ast$}}
\renewcommand{\@biblabel}[1]{#1.\hfill}
\makeatother

\renewcommand*{\leq}{\leqslant} \renewcommand*{\geq}{\geqslant}

\newcommand*{\mC}{\mathbb C} \newcommand*{\mQ}{\mathbb Q}
\newcommand*{\mR}{\mathbb R} \newcommand*{\mT}{\mathbb T}
\newcommand*{\mZ}{\mathbb Z}

\newcommand*{\cG}{\mathcal G} \newcommand*{\cO}{\mathcal O}
\newcommand*{\cT}{\mathcal T}

\newcommand*{\fE}{\mathfrak E} \newcommand*{\fT}{\mathfrak T}

\newcommand*{\barz}{\bar{z}}

\newcommand*{\gl}{\mathfrak{gl}} \newcommand*{\slsl}{\mathfrak{sl}}

\newcommand*{\GL}{\mathrm{GL}}

\newcommand*{\Georg}{\Gamma^{\#}} \newcommand*{\Cantor}{\Xi^{\#}}

\newcommand*{\bfomega}{\text{\boldmath $\omega^2$}}
\newcommand*{\bfsig}{\text{\boldmath $\sigma^n$}}
\newcommand*{\bfsigma}{\text{\boldmath $\sigma^{n+p}$}}

\newcommand*{\codim}{\mathop{\mathrm{codim}}\nolimits}

\newcommand*{\Ad}{\mathop{\mathrm{Ad}}\nolimits}
\newcommand*{\Fix}{\mathop{\mathrm{Fix}}\nolimits}
\newcommand*{\Spec}{\mathop{\mathrm{Spec}}\nolimits}
\newcommand*{\Tr}{\mathop{\mathrm{Tr}}\nolimits}

\begin{document}
\baselineskip=21pt
\parindent=0mm

\begin{center}
{\Large\bfseries Families of Invariant Tori in KAM Theory: \\
Interplay of Integer Characteristics}

\bigskip
\baselineskip=16pt

{\large\bfseries Mikhail B.~Sevryuk\footnote{E-mails: \texttt{2421584@mail.ru, sevryuk@mccme.ru}}}

\medskip

{\small\itshape V.~L.~Talroze Institute of Energy Problems of Chemical Physics of the Russia Academy of Sciences, \\
Leninskii prospect~38, Building~2, Moscow 119334, Russia}
\end{center}

\bigskip
\baselineskip=16pt

\textbf{Abstract}---%
The purpose of this brief note is twofold. First, we summarize in a very concise form the principal information on Whitney smooth families of quasi-periodic invariant tori in various contexts of KAM theory. Our second goal is to attract (via an informal discussion and a simple example) the experts' attention to the peculiarities of the so-called excitation of elliptic normal modes in the reversible context~2.

\bigskip

MSC2010 numbers: \texttt{37J40, 70H08, 70H33, 70K43}

\bigskip


\bigskip

\begin{flushright}\itshape
To the fond memory of Vladimir Igorevich Arnold, \\
one of the creators of contemporary mathematics
\end{flushright}

\parindent=5mm

\section{Whitney Smooth Families of Invariant Tori}
\label{Whitney}

Kolmogorov--Arnold--Moser (KAM) theory founded by the great mathematicians A.~N.~Kolmogorov (1903--1987), V.~I.~Arnold (1937--2010), and J.~K.~Moser (1928--1999) studies quasi-periodic motions in nonintegrable dynamical systems. The contribution of each of the three authors was described by Arnold in \cite{A14}. (Let me remark in parentheses that I was blessed by God enough to be Arnold's student at the Moscow State University in 1980--1987, to write my term papers, master's thesis, and PhD thesis under his supervision, and to learn KAM theory and the theory of reversible systems from him, as well as to discuss some aspects of KAM theory with Moser.) In this note, we will confine ourselves with KAM theory for autonomous dynamical systems with continuous time (i.e., autonomous flows) on finite dimensional manifolds. The central object of KAM theory is an invariant $n$-torus $\cT$ (to be more precise, an invariant manifold diffeomorphic to the standard $n$-torus $\mT^n=(\mR/2\pi\mZ)^n$) carrying \emph{quasi-periodic motions}. This means that in suitable coordinates $(x_1,\ldots,x_n)\in\mT^n$ in $\cT$, the equations of motion on $\cT$ take the ``linear'' form $\dot{x}=\omega$ with a constant vector $\omega\in\mR^n$, the components $\omega_1,\ldots,\omega_n$ of this vector being incommensurable (rationally independent). One also speaks of quasi-periodic invariant tori or nonresonant invariant tori. A more general concept is \emph{conditionally periodic motions} $\dot{x}=\omega$ for which the numbers $\omega_1,\ldots,\omega_n$ are arbitrary (not necessarily incommensurable). Invariant tori carrying conditionally periodic motions are also known as invariant tori with parallel dynamics or invariant tori with a Kronecker flow.

The main informal conclusion of KAM theory is that for many various classes of nonintegrable dynamical systems, quasi-periodic invariant $n$-tori with $n\geq 2$ are as typical and ``ubiquitous'' as invariant $0$-tori (equilibria) and invariant $1$-tori that do not contain equilibria (periodic trajectories, or cycles). The importance of quasi-periodic invariant tori (and, more generally, of invariant tori with parallel dynamics) stems, in the long run, from the fact that any finite dimensional connected and compact Abelian Lie group is a torus. Recent general reviews of KAM theory are exemplified by the tutorial \cite{dL01}, the monograph \cite[\S\S~6.2.2.C, 6.3]{AKN06}, and the survey \cite{BS10}. The book \cite{D14} presents a brilliant semi-popular introduction to the theory. The genericity approach to KAM theory with external parameters (this note is based on) has been developed mainly in the works \cite{BS10,BHT90,BH95,S95,BHS96Gro,BHS96LNM,S06,S07,BHN07,S16,S17} of the Groningen school of dynamical systems.

Quasi-periodic invariant tori (as well as arbitrary invariant tori with parallel dynamics) can be \emph{reducible} or nonreducible. Reducibility means that the variational equation along the torus can be reduced to a form with constant coefficients. Equilibria and periodic trajectories are always reducible.

Let $M$ be a smooth connected finite dimensional manifold (the phase space). We will consider a family of smooth vector fields $V_\mu$ on $M$ smoothly dependent on an \emph{external parameter} $\mu\in P\subset\mR^s$ ($s\geq 0$), where $P$ is an open domain. Generically, quasi-periodic invariant $n$-tori in KAM theory are organized into families which are smooth for $n=0,\,1$ and are \emph{Whitney smooth Cantor-like} for $n\geq 2$. The main ``ingredients'' of such a family of tori (in addition to $M$, $P$, and the family of vector fields $V_\mu$) in the reducible case are \cite{S95,BHS96Gro,BHS96LNM}:

a) an open domain $\Xi\subset\mR^c$ ($c\geq 0$) in which an \emph{internal parameter} $\nu$ ranges and a subset $\Cantor\subset\Xi$. For $n=0,\,1$ one has $\Cantor=\Xi$, whereas for $n\geq 2$ one always has $c\geq 1$ and the set $\Cantor$ is \emph{Cantor-like}, i.e., nowhere dense and of positive Lebesgue measure;

b) a smooth mapping
\begin{equation}
\Phi:\mT^n\times\cO_m(0)\times\Xi\to M\times P,
\label{Phi}
\end{equation}
where $m=\dim M-n\geq 0$ is the phase space codimension of the tori and $\cO_m(0)\subset\mR^m$ is a neighborhood of the origin in $\mR^m$. This mapping possesses the following properties (below $x$, $X$, and $\nu$ are coordinates in $\mT^n$, $\cO_m(0)$, and $\Xi$, respectively, while $w$ is a point in $M$).

First, the restriction of $\Phi$ to $\mT^n\times\{0\}\times\Cantor$ is injective.

Second, for any $\nu\in\Xi$, the set $\Phi\bigl( \mT^n\times\cO_m(0)\times\{\nu\} \bigr)$ lies in one of the fibers $M\times\{\mu^\nu\}$ of the bundle $M\times P\to P$, $(w,\mu)\mapsto\mu$. Thus, the mappings
\[
\Phi^\nu:\mT^n\times\cO_m(0)\to M, \quad \Phi(x,X,\nu)=\bigl( \Phi^\nu(x,X), \, \mu^\nu \bigr)
\]
are well-defined.

Third, for any $\nu\in\Xi$, the mapping $\Phi^\nu:\mT^n\times\cO_m(0)\to M$ is a diffeomorphism onto its image.

Fourth, for any $\nu\in\Cantor$, the vector field $V_{\mu^\nu}$ in the coordinates $(x,X)$ induced by $\Phi^\nu$ (the so-called \emph{Floquet coordinates}) affords the equations of motion
\begin{equation}
\dot{x}=\omega^\nu+O(X), \quad \dot{X}=\Lambda^\nu X+O\bigl( |X|^2 \bigr),
\label{Floq}
\end{equation}
where $\omega^\nu\in\mR^n$ is a certain constant vector and $\Lambda^\nu\in\gl(m,\mR)$ is a certain constant matrix.

Fifth, for all $\nu\in\Cantor$, the components $\omega^\nu_1,\ldots,\omega^\nu_n$ of the vector $\omega^\nu$ are uniformly strongly incommensurable, i.e., are independent over rationals and uniformly badly approximable by rationally dependent quantities. One may regard $\omega^\nu_1,\ldots,\omega^\nu_n$ as being uniformly \emph{Diophantine}. This means the existence of constants $\tau>\max(n-1,0)$ and $\gamma>0$ independent of $\nu\in\Cantor$ and such that
\begin{equation}
\bigl| \langle\omega^\nu,k\rangle \bigr| \geq \gamma|k|^{-\tau} \quad \forall\; k\in\mZ^n\setminus\{0\}
\label{Dioph}
\end{equation}
(here and henceforth, the angle brackets denote the standard inner product of vectors).

For each $\nu\in\Cantor$, the set $\cT^\nu=\Phi^\nu\bigl( \mT^n\times\{0\} \bigr)=\{X=0\}$ is a reducible invariant $n$-torus of the flow of the vector field $V_{\mu^\nu}$. The motion on this torus is quasi-periodic with strongly incommensurable (say, Diophantine) frequencies $\omega^\nu_1,\ldots,\omega^\nu_n$. Besides the frequency vector $\omega^\nu$, the torus $\cT^\nu$ is characterized by the \emph{Floquet matrix} $\Lambda^\nu$. The eigenvalues of this matrix are called the \emph{Floquet exponents} of the torus $\cT^\nu$, and the positive imaginary parts $\beta^\nu_1,\ldots,\beta^\nu_q$ of the Floquet exponents ($q$ lies in the interval $0\leq q\leq m/2$ and is independent of $\nu$) are called the \emph{normal frequencies} of the torus. Correspondingly, the numbers $\omega^\nu_1,\ldots,\omega^\nu_n$ are sometimes called the \emph{tangential} frequencies of the torus $\cT^\nu$. As a rule, one deals with families of tori whose tangential and normal frequencies are uniformly \emph{affinely Diophantine} with bound $2$. This means that for some constants $\tau>\max(n-1,0)$ and $\gamma>0$ independent of $\nu\in\Cantor$, one has
\[
\bigl| \langle\omega^\nu,k\rangle+\langle\beta^\nu,l\rangle \bigr| \geq \gamma|k|^{-\tau} \quad \forall\; k\in\mZ^n\setminus\{0\}, \; l\in\mZ^q, \; |l|=|l_1|+\cdots+|l_q|\leq 2.
\]

For $\nu\in\Xi\setminus\Cantor$, the mapping $\Phi^\nu$ and the coordinates $(x,X)$ in $M$ induced by $\Phi^\nu$ have no dynamical meaning. The motion in the corresponding ``gaps'' between the tori $\cT^\nu$ is usually very complicated, and these gaps often contain other Whitney smooth families of invariant tori of various dimensions (cf.\ Section~\ref{destruction}). The factor $\cO_m(0)$ in~\eqref{Phi} is absent in the nonreducible case.

Small linear combinations $\langle\omega^\nu,k\rangle$ of the frequencies $\omega^\nu_1,\ldots,\omega^\nu_n$ in~\eqref{Dioph} for large $k\in\mZ^n$ are the famous \emph{small divisors} which are the main source of difficulties in the theory of quasi-periodic motions. Of course, small divisors are impossible for $n\leq 1$, and that is why the cases $n=0,\,1$ and $n\geq 2$ in KAM theory are so different. In fact, equilibria and periodic trajectories are often regarded as being outside the scope of KAM theory.

\section{Contexts of KAM Theory}
\label{contexts}

The particular structure of generic families $\bigl\{ \cT^\nu \bigm| \nu\in\Cantor \bigr\}$ of quasi-periodic invariant tori depends on the conservation laws the vector fields $V_\mu$ are assumed to obey and on the symmetry properties these vector fields are assumed to possess. Such conservation laws and symmetry properties constitute what is called the \emph{context} of KAM theory. The four best explored KAM contexts are the following ones \cite{BS10,BHT90,BH95,S95,BHS96Gro,BHS96LNM,S06,S07,BHN07,S16,S17,W10} ($n\geq 0$ always denotes the dimension of the quasi-periodic invariant tori under consideration, while $s$ is the number of external parameters $\mu_1,\ldots,\mu_s$).

1) The \emph{Hamiltonian isotropic} context, where the phase space $M$ is a symplectic manifold equipped with a symplectic structure $\bfomega$, the vector fields $V_\mu$ are Hamiltonian (in other words, the interior products $i_{V_\mu}\bfomega$ of $\bfomega$ with $V_\mu$ are exact: $i_{V_\mu}\bfomega=dH_\mu$), and the invariant tori $\cT^\nu$ are \emph{isotropic}, i.e., the restrictions of $\bfomega$ to $\cT^\nu$ vanish. The isotropy of the $n$-tori $\cT^\nu$ implies that $2n\leq\dim M$. If $\dim M=2(n+p)$ with $p\geq 0$, we will speak of the Hamiltonian isotropic $(n,p,s)$ context.

2) The \emph{volume preserving} $(n,p,s)$ context ($p\geq 1$), where the phase space $M$ of dimension $n+p$ is equipped with a volume element $\bfsigma$ and the vector fields $V_\mu$ are globally volume preserving, or globally divergence-free, i.e., the $(n+p-1)$-forms $i_{V_\mu}\bfsigma$ are exact. The volume preserving $(n,0,s)$ context is impossible \cite{BHS96Gro,BHS96LNM}. This can be easily explained heuristically as follows. Consider the torus $\mT^n\ni x$ equipped with the volume element $\bfsig=dx_1\wedge\cdots\wedge dx_n$. Then any nonzero constant vector field $V=\omega\partial/\partial x$ on $\mT^n$ is volume preserving ($di_V\bfsig=0$) but not \emph{globally} volume preserving. In fact, the correspondence ``$\omega$ $\mapsto$ the cohomology class of $i_V\bfsig$'' determines an isomorphism $\mR^n\to H^{n-1}(\mT^n,\mR)$.

3) The \emph{general dissipative} $(n,p,s)$ context ($p\geq 0$), where the phase space $M$ of dimension $n+p$ is not assumed to be equipped with any special structure.

4) The \emph{reversible} context, where a smooth involution $G:M\to M$ (a mapping whose square is the identical transformation) of the phase space $M$ is given, the vector fields $V_\mu$ are reversible with respect to $G$, and the tori $\cT^\nu$ invariant under the flows of $V_{\mu^\nu}$ are also invariant under $G$. The reversibility of $V_\mu$ with respect to $G$ means that the involution $G$ casts the fields $V_\mu$ into the opposite fields $-V_\mu$:
\[
\Ad_GV_\mu=TG(V_\mu\circ G^{-1})=TG(V_\mu\circ G)=-V_\mu.
\]
It is well known that the fixed point set
\[
\Fix G=\bigl\{ w\in M \bigm| G(w)=w \bigr\}
\]
of any smooth involution $G:M\to M$ is a submanifold of $M$ of the same smoothness class as the involution $G$ itself. Moreover, if a quasi-periodic invariant $n$-torus $\cT$ of a $G$-reversible vector field $V$ is also invariant under $G$, then one can choose a coordinate $x\in\mT^n$ in $\cT$ in such a way that the dynamics on $\cT$ will take the form $\dot{x}=\omega$ and the restriction of $G$ to $\cT$ will have the form $G|_{\cT}:x\mapsto-x$ \cite{BHS96Gro,BHS96LNM}. Consequently, the set $(\Fix G)\cap\cT=\Fix\bigl( G|_{\cT} \bigr)$ consists of $2^n$ isolated points
\[
(x_1,\ldots,x_n), \quad x_j\in\{0;\pi\}, \; 1\leq j\leq n,
\]
and the codimension of any connected component of $\Fix G$ that intersects $\cT$ is no less than $n$. In the reversible context of KAM theory, one usually assumes that all the connected components of $\Fix G$ are of the same dimension (the case most often encountered in practice), so that the numbers $\dim\Fix G=a\geq 0$ and $\codim\Fix G=n+b$ ($b\geq 0$) are well-defined. In such a setup with $\dim M=a+b+n$, we will speak of the reversible $(n,a,b,s)$ context. One usually distinguishes \emph{the reversible context~1} where
\[
a\geq b \Longleftrightarrow 2a=2\dim\Fix G \geq \codim\cT=a+b
\]
and the \emph{the reversible context~2} where
\[
a<b \Longleftrightarrow 2a=2\dim\Fix G < \codim\cT=a+b,
\]
see \cite{BHS96Gro,BHS96LNM,S16,S17} and references in \cite{S16,S17}.

Within the reversible $(n,a,b,s)$ context, we also require in the construction of Section~\ref{Whitney} that for any $\nu\in\Xi$, the reversing involution $G$ in the coordinates $(x,X)$ induced by $\Phi^\nu$ has the form $G:(x,X)\mapsto(-x,RX)$ where $R\in\GL(a+b,\mR)$ is an involutive matrix independent of $\nu$. The eigenvalues $1$ and $-1$ of $R$ are of multiplicities $a$ and $b$, respectively.

If an invariant torus $\cT$ of a $G$-reversible vector field $V$ is not invariant under $G$ then $G(\cT)\neq\cT$ is also invariant under the flow of $V$ but the dynamics near $\cT$ and $G(\cT)$ is essentially dissipative. In fact, coexistence of regions in the phase space with conservative-like dynamics and those with dissipative-like dynamics is a typical phenomenon for reversible systems \cite{POB86,QR89}.

The smoothness class of families of invariant tori is determined by that of the families of vector fields. For instance, consider the analytic category where the phase space $M$, the structures $\bfomega$, $\bfsigma$, or $G$ on $M$, the vector fields $V_\mu$, and their dependence on the external parameter $\mu$ are real analytic. Then the invariant $n$-tori $\cT^\nu$ themselves are also analytic whereas the family they constitute is analytic for $n=0,\,1$ and is $C^\infty$-smooth \emph{in the sense of Whitney} for $n\geq 2$. To be more precise, the mapping $\Phi$~\eqref{Phi} is analytic in $x\in\mT^n$ and $X\in\cO_m(0)$ for any $n\geq 0$, is analytic in $\nu\in\Xi$ for $n=0,\,1$, and is infinitely differentiable in $\nu\in\Xi$ for $n\geq 2$. In fact, $\Phi$ is Gevrey regular in $\nu\in\Xi$ for $n\geq 2$, see the papers \cite{W10,W03,S03} and references therein. In the $C^\infty$-category, the mapping $\Phi$~\eqref{Phi} is of class $C^\infty$ in all its arguments.

For basic references on Whitney smoothness in KAM theory, see \cite{BS10,BHS96LNM}. Whitney smoothness of families of invariant tori in the reversible context~2 was proven only in 2016 \cite{S16}.

Some essential features of generic smooth (for $n=0,\,1$) and Whitney smooth Cantor-like (for $n\geq 2$) families of quasi-periodic invariant $n$-tori in KAM theory \cite{BS10,S95,BHS96Gro,BHS96LNM,S06,S07,S17} are summarized in Tables~1 and~2. To be more precise, in the functional space of all the families $V_\mu$ of vector fields within the context in question, there is an open subset of families of vector fields admitting quasi-periodic motions with principal characteristics indicated in Tables~1 and~2. Families of invariant tori with other properties can only be encountered as an exception. In other words, such families either indicate the existence of some additional (explicit or hidden) symmetries in the systems or exhibit smaller values of the number $c$ of internal parameters (and are ``adjacent'' to generic families of invariant tori in the same sense as a complicated singularity can be adjacent to a simpler singularity). For instance, the paper \cite{BCHV09} is devoted to reducible invariant tori in the reversible $(n,a,b,s)$ context~1 where the multiplicity of zero Floquet exponent of the tori is larger than $a-b$ but the number of internal parameters of the Whitney smooth family of the tori is smaller than $a-b+s$ (see also \cite{S16} for a discussion). The Floquet matrices of the invariant tori in \cite{BCHV09} are not diagonalizable over $\mC$: their Jordan structure involves at least one nilpotent Jordan block of order greater than $1$.

\begin{table}[tb]
\begin{center}\small\renewcommand{\arraystretch}{1.5}
\textbf{Table~1.} Key properties of generic families of invariant $n$-tori $\cT$ in various KAM contexts. The numbers $n$, $p$, $a$, and $b$ can be any nonnegative integers unless stated otherwise. The meaning of the integers $s$, $p$, $a$, $b$, $m$, and $c$ is explained in the text, and $M$ is the phase space. The number $\delta_{1p}$ is equal to $1$ for $p=1$ and to $0$ for $p\geq 2$. The number $\Delta_n$ is equal to $0$ for $n=0,\,1$ and to $1$ for $n\geq 2$.

\medskip

\begin{tabular}{|c|c|c|c|c|c|}\hline
Context & $\dim M$ & $m$ & Lower bound of $s$ & $c$ & $c-s$ \\ \hline
Hamiltonian isotropic $(n,p,s)$ & $2(n+p)$ & $n+2p$ & $s\geq 0$ & $n+s$ & $n$ \\ \hline
Volume preserving $(n,p,s)$, $p\geq 1$ & $n+p$ & $p$ & $s\geq\max(\Delta_n-\delta_{1p},0)$ & $\delta_{1p}+s$ & $\delta_{1p}$ \\ \hline
General dissipative $(n,p,s)$ & $n+p$ & $p$ & $s\geq\Delta_n$ & $s$ & $0$ \\ \hline
Reversible $(n,a,b,s)$ & $n+a+b$ & $a+b$ & $s\geq\max(\Delta_n+b-a,0)$ & $a-b+s$ & $a-b$ \\ \hline
\end{tabular}
\end{center}
\end{table}

\begin{table}[tb]
\begin{center}\small\renewcommand{\arraystretch}{1.5}
\textbf{Table~2.} The spectra of the $m\times m$ Floquet matrices $\Lambda$ of reducible invariant $n$-tori $\cT$ in various KAM contexts. Generically each matrix $\Lambda$ has $g$ nonzero eigenvalues, and the multiplicity of the zero eigenvalue is equal to $|c-s|$.

\medskip

\begin{tabular}{|c|c|c|}\hline
Context & $g=m-|c-s|$ & Floquet matrices $\Lambda$ \\ \hline
Hamiltonian isotropic $(n,p,s)$ & $2p$ & $\Spec\Lambda = \bigl\{ \,\underbrace{0,\ldots,0}_n\,, \; \pm\lambda_1,\ldots,\pm\lambda_p \bigr\}$ \\ \hline
Volume preserving $(n,p,s)$, $p\geq 1$ & $p-\delta_{1p}$ & $\Tr\Lambda=0$ \\ \hline
General dissipative $(n,p,s)$ & $p$ & nothing special \\ \hline
Reversible $(n,a,b,s)$ & $2\min(a,b)$ & $\Spec\Lambda = \bigl\{ \,\underbrace{0,\ldots,0}_{|a-b|}\,, \; \pm\lambda_1,\ldots,\pm\lambda_{\min(a,b)} \bigr\}$ \\ \hline
\end{tabular}
\end{center}
\end{table}

Recall that in Tables~1 and~2, $s$ is the number of external parameters $\mu_1,\ldots,\mu_s$ (the dimension of the domain $P$ where the external parameter $\mu$ ranges) while $c$ is the number of internal parameters $\nu_1,\ldots,\nu_c$ (the dimension of the domain $\Xi$ where the internal parameter $\nu$ ranges). One speaks of $s$-parameter families of vector fields $V_\mu$ and $c$-parameter families of quasi-periodic invariant tori $\cT^\nu$. The really essential information in Tables~1 and~2 is the lower bound of $s$ (the minimal value of $s$ for which quasi-periodic invariant $n$-tori $\cT$ occur generically), the value of $c$, and the properties of the spectra of the Floquet matrices $\Lambda$ of reducible tori $\cT$.

The letter $m$ in Tables~1 and~2 always denotes the phase space codimension $\dim M-n$ of the invariant $n$-tori $\cT^\nu$ in question. The symbol $\delta_{ij}$ is the Kronecker delta, that is, $\delta_{ij}=1$ for $i=j$ and $\delta_{ij}=0$ for $i\neq j$. The symbol $\Delta_n$ is $1-\delta_{0n}-\delta_{1n}$, i.e., $\Delta_n=0$ for $n=0,\,1$ and $\Delta_n=1$ for $n\geq 2$.

The lower bound of $s$ in Table~1 is determined by the condition $c\geq\Delta_n$ in all the four contexts. Quasi-periodic invariant $n$-tori for $n\geq 2$ are not isolated in the product $M\times P$ of the phase space and the external parameter space, such tori are generically organized into at least one-parameter Whitney smooth families. In all the cases presented in Tables~1 and~2 except for the reversible context~2, the values $s=0$ or $s=1$ are enough for the generic existence of families of invariant tori in question \cite{S17}.

The meaning of the last column of Table~1 is as follows. If $c>s$ then $c-s$ is the number of parameters of the family of the tori $\cT^\nu$ corresponding to each \emph{individual} vector field $V_{\mu_0}$ (i.e., of the tori $\cT^\nu$ with $\mu^\nu=\mu_0$). In other words, $c-s\leq m$ and each individual vector field $V_\mu$ (for $\mu$ in some open subdomain of $P$) admits generically a $(c-s)$-parameter family of quasi-periodic invariant $n$-tori (according to the general rule, this family is smooth for $n=0,\,1$ and is Whitney smooth Cantor-like for $n\geq 2$). If $c\leq s$ then the map
\[
\mR^c\supset\Xi\ni\nu\mapsto\mu^\nu\in P\subset\mR^s
\]
is generically a (local) diffeomorphism onto its image $\Gamma\subset P$. The set $\Gamma$ is a smooth surface of codimension $s-c$ in $P$. Let
\[
\Georg=\bigl\{ \mu^\nu \bigm| \nu\in\Cantor \bigr\}\subset\Gamma.
\]
If $\mu\in\Georg$ then the vector field $V_\mu$ possesses a \emph{single} invariant $n$-torus of the given family $\cT^\nu$. On the other hand, if $\mu\notin\Georg$ then the vector field $V_\mu$ admits \emph{no} invariant $n$-tori of the given family.

If $n\geq 2$, $c>s$, and $c-s=m$, then the ``gaps'' between the invariant tori $\cT^\nu$ in the phase space are called resonant zones.

According to the last column of Table~2, generic reducible invariant tori $\cT^\nu$ in the Hamiltonian isotropic $(n,p,s)$ context have the zero Floquet exponent of multiplicity $n=c-s$ (for $n\geq 1$) while the remaining $2p$ Floquet exponents come in pairs $\pm\lambda$. Similarly, generic reducible invariant tori $\cT^\nu$ in the reversible $(n,a,b,s)$ context have the zero Floquet exponent of multiplicity $|a-b|=|c-s|$ (for $a\neq b$) while the remaining $2\min(a,b)$ Floquet exponents come in pairs $\pm\lambda$. The traces of the $p\times p$ Floquet matrices $\Lambda^\nu$ of reducible invariant tori $\cT^\nu$ in the volume preserving $(n,p,s)$ context vanish. For $p=1$ (where $c-s=1$) this implies $\Lambda^\nu\equiv 0$. On the other hand, generic matrices in $\slsl(p,\mR)$ for $p\geq 2$ (where $c=s$) are nonsingular. Of course, generic matrices in $\gl(p,\mR)$ for any $p$ are also nonsingular. One sees that the multiplicity of the zero Floquet exponent of reducible invariant tori is generically equal to $|c-s|$ in all the four KAM contexts, and the number of nonzero Floquet exponents is generically equal to $g=m-|c-s|$.

Note that in the volume preserving $(n,p,s)$ context, $c-s=\delta_{1p}$ is equal to the $(n+p-1)$-th Betti number of $\mT^n\times\mR^p$ for each $p\geq 1$.

Invariant tori $\cT^\nu$ in the Hamiltonian isotropic $(n,p,s)$ context with $p=g/2\geq 1$ and in the reversible $(n,a,b,s)$ context with $\min(a,b)=g/2\geq 1$ are often said to be \emph{lower-dimensional}.

Of course, for the Hamiltonian isotropic context and the reversible context~1, the information compiled in Tables~1 and~2 is in fact rather well known.

\section{Proper Destruction of Resonant Tori}
\label{destruction}

It is amazing that the very limited and formal data presented in Tables~1 and~2 enable one to predict, by means of just ``juggling with integers'', essential features of some quite complicated phenomena studied in KAM theory.

One of such phenomena is destruction of resonant tori. Let $n\geq 2$ and let $r$ be an integer in the range $1\leq r\leq n-1$. Within any of the KAM contexts of Section~\ref{contexts}, consider partially integrable vector fields $V_\mu$ for which the corresponding family of invariant $n$-tori $\cT^\nu$ is smooth rather than Whitney smooth Cantor-like. The tori $\cT^\nu$ carry conditionally periodic motions with frequency vectors $\omega^\nu$. Suppose that $c\geq r$. Then, typically, the smooth $c$-parameter family of $n$-tori $\cT^\nu$ contains a smooth $(c-r)$-parameter subfamily of tori whose frequencies satisfy $r$ independent \emph{fixed} resonance relations
\[
\Bigl\langle \omega^\nu,k^{(\iota)} \Bigr\rangle=0, \quad 1\leq\iota\leq r,
\]
where $k^{(1)},\ldots,k^{(r)}\in\mZ^n\setminus\{0\}$ are fixed integer vectors linearly independent over $\mQ$. Such a resonance is usually said to be simple for $r=1$ and multiple for $r>1$. Each resonant $n$-torus is foliated into invariant $(n-r)$-tori with parallel dynamics. We will say that in this setup, \emph{proper destruction of resonant tori} takes place if, under a generic small perturbation of $V_\mu$ within the given context, this smooth $(c-r)$-parameter family of resonant invariant $n$-tori gives rise to a \emph{finite collection} of $(c-r)$-parameter families of quasi-periodic invariant $(n-r)$-tori, the latter families being smooth for $r=n-1$ and Whitney smooth Cantor-like for $r<n-1$ and $c>r$. The number of such families and the normal behavior of the perturbed invariant $(n-r)$-tori depend on the perturbation.

In the terms of the data of Table~1, proper destruction of resonant tori means a passage from $n$ to $n'=n-r$ in such a way that the values of $\dim M$ and $s$ (as well as the value of $a=\dim\Fix G$ in the reversible context) remain unchanged while the value of $c\geq r$ decreases by $r$ (together with $n$). Table~1 allows one to conjecture in what cases proper destruction of resonant tori is possible.

In the Hamiltonian isotropic $(n,p,s)$ context, proper destruction of resonant tori can occur for any $n\geq 2$, any $r$ in the range $1\leq r\leq n-1$, and any nonnegative $p$ and $s$. Indeed, a passage from $n$ and $p$ to $n'=n-r$ and $p'=p+r$ preserves $\dim M=2(n+p)$ and yields $c'=n'+s=c-r$. Moreover, one always has $c'\geq 1$. Proper destruction of resonant tori in Hamiltonian systems has been really described and explored, but (to the best of the author's knowledge) for $p=0$ and $s=0$ only. For $r=n-1$ (i.e., for $n'=1$) the corresponding perturbed cycles are sometimes called \emph{Poincar\'e trajectories} because the research into them goes back to H.~Poincar\'e, see \cite[\S~4.2.1]{BHS96LNM} and references therein. The studies of proper destruction of resonant tori in Hamiltonian systems for $r\leq n-2$ (i.e., for $n'\geq 2$) were started by D.~V.~Treshch\"ev \cite{T89}, so that the corresponding perturbed $(n-r)$-tori are called \emph{Treshch\"ev tori} \cite[\S~4.2.2]{BHS96LNM} or Poincar\'e--Treshch\"ev tori \cite{LY05}. In the case of a simple resonance ($r=1$ and $n'=n-1$), the terms ``Birkhoff--Kolmogorov--Arnold--Moser tori'' \cite{C96} or ``Birkhoff lower-dimensional tori'' \cite{WC97} are also used. For other references on proper destruction of resonant tori in the Hamiltonian isotropic context and a survey of the main results, see \cite[p.~305]{AKN06}, \cite[Section~8.2]{BS10}, and \cite{S03}.

In the volume preserving $(n,p,s)$ context, the value of $c=\delta_{1p}+s$ can decrease for $s$ fixed only if one passes from $p=1$ to $p'>1$. In this case $c'=s=c-1$ and $\dim M=n+p=n+1$ whence $r=c-c'=1$, $n'=n-r=n-1$, and $p'=\dim M-n'=2$. Thus, in this context, proper destruction of resonant tori can happen for any $n\geq 2$, but only with $p=1$ and $r=1$. In this case, as we have just seen, $c=s+1$, $n'=n-1$, $p'=2$, and $c'=s$. The condition $c'\geq\Delta_{n'}$ implies that for $n\geq 3$ (i.e., for $n'\geq 2$) one should additionally require $s\geq 1$. However, it seems that proper destruction of ``simply resonant'' tori of codimension $1$ in volume preserving flows has not been examined yet.

In the general dissipative context, proper destruction of resonant tori is impossible because here $c=s$ cannot change for $s$ fixed.

In the reversible $(n,a,b,s)$ context, proper destruction of resonant tori can take place for any $n\geq 2$, any $r$ in the range $1\leq r\leq n-1$, any nonnegative $a$ and $b$, and any $s\geq\max(\Delta_{n-r}+b+r-a,0)$. Indeed, a passage from $n$ and $b$ to $n'=n-r$ and $b'=b+r$ with $a$ unaltered preserves $\dim M=n+a+b$ and yields $c'=a-b'+s=a-b-r+s=c-r$. The condition $c'\geq\Delta_{n'}$ is tantamount to the inequality $s\geq\Delta_{n-r}+b+r-a$. As far as the author knows, by now proper destruction of resonant tori in reversible systems has been described and analyzed for $a=n$, $b=0$, and $s=0$ only \cite{L01,W01}. In this case $c=n$, $b'=r$, $c'=n-r=n'$, one always has $c'\geq 1$, the unperturbed invariant $n$-tori pertain to the reversible context~1, and so do the perturbed invariant $n'$-tori.

\section{Excitation of Elliptic Normal Modes}
\label{excitation}

Another phenomenon Tables~1 and~2 very much help to understand is excitation of so-called elliptic normal modes. According to Table~2, in the Hamiltonian isotropic $(n,p,s)$ context with $p\geq 1$, in the volume preserving $(n,2,s)$ context, and in the reversible $(n,a,b,s)$ context with $\min(a,b)\geq 1$, reducible invariant $n$-tori $\cT^\nu$ generically have nonzero Floquet exponents, and these exponents come in pairs $\pm\lambda$. Consequently, in the contexts indicated, some Floquet exponents of reducible invariant $n$-tori $\cT^\nu$ can be nonzero purely imaginary with ``a positive probability''. This is no longer so in the volume preserving $(n,p,s)$ context for $p\geq 3$ or in the general dissipative context. Indeed, if $\Lambda\in\slsl(2,\mR)$ then $\Spec\Lambda=\{\pm\lambda\}$, but generic matrices in $\slsl(p,\mR)$ for $p\geq 3$ have no pairs $\pm\lambda$ of opposite eigenvalues.

Now, within the framework of the Hamiltonian isotropic $(n,p,s)$ context with $p\geq 1$, the volume preserving $(n,2,s)$ context, or the reversible $(n,a,b,s)$ context with $\min(a,b)\geq 1$, consider partially integrable vector fields $V_\mu$ admitting a smooth $c$-parameter family of reducible invariant $n$-tori $\cT^\nu$ with parallel dynamics. Let $\omega^\nu_1,\ldots,\omega^\nu_n$ be the frequencies of these tori. Suppose that among the $g$ nonzero Floquet exponents of each of the tori $\cT^\nu$ (see Table~2), there are $r\geq 1$ pairs of purely imaginary numbers $\pm i\beta^\nu_1,\ldots,\pm i\beta^\nu_r$, where $r\leq p$ in the Hamiltonian isotropic $(n,p,s)$ context, $r=1$ in the volume preserving $(n,2,s)$ context, and $r\leq\min(a,b)$ in the reversible $(n,a,b,s)$ context. The remaining $g-2r$ nonzero Floquet exponents of $\cT^\nu$ are also allowed to be purely imaginary. One says that the Floquet exponents $\pm i\beta^\nu_1,\ldots,\pm i\beta^\nu_r$ ``excite'' if, near the family of the $n$-tori $\cT^\nu$, the vector fields $V_\mu$ themselves and any sufficiently small perturbations of $V_\mu$ within the given context possess a $(c+r)$-parameter family of reducible quasi-periodic invariant $(n+r)$-tori $\fT$ with frequencies close to
\[
\omega^\nu_1,\ldots,\omega^\nu_n, \beta^\nu_1,\ldots,\beta^\nu_r.
\]
The latter family is smooth for $n=0$, $r=1$ and is Whitney smooth Cantor-like for $n+r\geq 2$. This is what is called \emph{excitation of elliptic normal modes}.

In a sense, excitation of elliptic normal modes is a phenomenon \emph{opposite} to proper destruction of resonant tori. In the terms of the data of Table~1, excitation of elliptic normal modes means a passage from $n$ to $n'=n+r$ in such a way that the values of $\dim M$ and $s$ (as well as the value of $a=\dim\Fix G$ in the reversible context) remain unchanged while the value of $c$ increases by $r$ (together with $n$). In the Hamiltonian isotropic $(n,p,s)$ context (where $r\leq p$), this is achieved by a passage from $p$ to $p'=p-r$. In the volume preserving $(n,2,s)$ context (where $r=1$), one passes from $p=2$ to $p'=1$. In the reversible $(n,a,b,s)$ context (where $r\leq\min(a,b)$), a passage from $b$ to $b'=b-r$ takes place. By the way, if the initial $n$-tori $\cT^\nu$ pertain to the reversible context~1 (i.e., if $a\geq b$), so do the $(n+r)$-tori $\fT$ a~fortiori.

To get the picture of excitation of elliptic normal modes, one may think of a very particular case where the vector fields $V_{\mu^\nu}$ afford the equations of motion
\begin{equation}
\begin{aligned}
\dot{x} &= \omega^\nu, \\
\dot{z}_j &= i\beta^\nu_jz_j\bigl[ 1+F^\nu_j(z_1\barz_1,\ldots,z_r\barz_r) \bigr], \quad j=1,\ldots,r, \\
\dot{y} &= L^\nu y
\end{aligned}
\label{Birk}
\end{equation}
near every $n$-torus $\cT^\nu=\{z=0, \: y=0\}$, where $x\in\mT^n$, the variable $z$ ranges in a neighborhood $\cO_r^{\mC}(0)$ of the origin in $\mC^r$, the variable $y$ ranges in a neighborhood $\cO_{m-2r}(0)$ of the origin in $\mR^{m-2r}$, the functions $F^\nu_j$ are real-valued, $F^\nu_j(0)=0$ for each $j$, and $L^\nu\in\gl(m-2r,\mR)$. For the Floquet matrices $\Lambda^\nu$ of the tori $\cT^\nu$, one has
\[
\Spec\Lambda^\nu = \{\pm i\beta^\nu_1;\ldots;\pm i\beta^\nu_r\}\cup\Spec L^\nu.
\]
In this situation, from every $n$-torus $\cT^\nu$, there ``emanates'' the $r$-parameter smooth family of invariant $(n+r)$-tori
\begin{equation}
\{z_1\barz_1=C_1>0, \: \ldots, \: z_r\barz_r=C_r>0, \: y=0\}
\label{C1Cr}
\end{equation}
with parallel dynamics and with the frequencies
\[
\omega^\nu_1,\ldots,\omega^\nu_n, \beta^\nu_1\bigl[ 1+F^\nu_1(C_1,\ldots,C_r) \bigr],\ldots,\beta^\nu_r\bigl[ 1+F^\nu_r(C_1,\ldots,C_r) \bigr].
\]
Under a perturbation, such a smooth family is expected to become Whitney smooth Cantor-like (except for the case where $n=0$ and $r=1$).

In the Hamiltonian isotropic $(n,p,s)$ context with $p\geq 1$, in the volume preserving $(n,2,s)$ context, and in the reversible $(n,a,b,s)$ context~1 with $a\geq b\geq 1$, excitation of elliptic normal modes has been described and studied in detail. For $n=0$, $r=1$, and $s=0$, smooth one-parameter families of cycles $\fT$ in the Hamiltonian isotropic $(0,p,0)$ context are called Lyapunov families, and those in the reversible $(0,a,a,0)$ context~1 are called Lyapunov--Devaney families, see e.g.\ \cite[p.~96]{BHS96LNM} and \cite{D76}. Excitation of elliptic normal modes for the nontrivial case $n'=n+r\geq 2$ was first observed by V.~I.~Arnold \cite{A61} (in the Hamiltonian isotropic $(0,2,0)$ context with $r=2$), see his recollections \cite{A14}. In fact, for $n=0,\,1$, the question is the dynamics in a neighborhood of an equilibrium or periodic trajectory with partially elliptic normal behavior (and the stability of such equilibria or periodic trajectories). This is the subject of the so-called local KAM theory dealt with in an enormous number of works.

Excitation of elliptic normal modes in the more complicated setting where $n\geq 2$ was also first handled by Arnold \cite{A62,A63}, mainly in connection with the problem of stability of planetary systems. In our notation, in \cite{A62,A63} Arnold considered excitation of elliptic normal modes in the Hamiltonian isotropic $(n,p,0)$ context for any $n$ and $p$ with $r=p$. The invariant $(n+r)$-tori $\fT$ in \cite{A61,A62,A63} were Lagrangian: they were isotropic and their dimension was equal to the number $n+p$ of degrees of freedom. Arnold's results \cite{A62,A63} were recently revised, extended, and strengthened in an essential way \cite{F04,CP09,CP10,CP11IM,F13,P15,Pnow} (see also \cite{CP11JMD,P13}). Excitation of elliptic normal modes in the Hamiltonian isotropic $(n,p,0)$ context for any $n$, $p$, and $r\leq p$ and in the reversible $(n,n+b,b,0)$ context~1 for any $n$, $b$, and $r\leq b$ was conjectured in \cite{S90}. All the references below will pertain to excitation of purely imaginary Floquet exponents of invariant $n$-tori with $n\geq 2$.

In all the three KAM contexts indicated above (the Hamiltonian isotropic context, the volume preserving context with $p=2$, and the reversible context~1), excitation of elliptic normal modes is characterized by the following features. First, the number of zero Floquet exponents of the $(n+r)$-tori $\fT$ is larger by $r$ than that of the initial $n$-tori $\cT^\nu$. Indeed, in all the three contexts in question, one has $c\geq s$ (see Table~1), so that the number of zero Floquet exponents of invariant tori is equal to $|c-s|=c-s$, and an increase in $c$ by $r$ implies an increase in $|c-s|$ by $r$. Second, if the nonzero Floquet exponents of the initial $n$-tori $\cT^\nu$ are
\[
\pm i\beta^\nu_1,\ldots,\pm i\beta^\nu_r, \pm\lambda^\nu_1,\ldots,\pm\lambda^\nu_{(g-2r)/2}, \quad g=m-c+s
\]
(see Table~2), then the nonzero Floquet exponents of the $(n+r)$-tori $\fT$ are close to
\[
\pm\lambda^\nu_1,\ldots,\pm\lambda^\nu_{(g-2r)/2}.
\]
Third, to each individual vector field, there generically corresponds a $(c+r-s)$-parameter family of invariant $(n+r)$-tori $\fT$. This can be understood as follows: each individual unperturbed vector field $V_\mu$ has a $(c-s)$-parameter smooth family of the initial $n$-tori $\cT^\nu$, and each such torus ``emits'' an $r$-parameter subfamily of invariant $(n+r)$-tori $\fT$.

Excitation of the Floquet exponents $\pm i\beta^\nu_1,\ldots,\pm i\beta^\nu_r$ of the initial invariant $n$-tori $\cT^\nu$ in all the three contexts indicated above is implied by certain nondegeneracy and nonresonance conditions imposed on the vector fields $V_\mu$. If $c\geq 1$, such nondegeneracy conditions can be formulated exclusively in terms of the frequencies and Floquet exponents of the tori $\cT^\nu$, see e.g.\ \cite[\S~4.1.4]{BHS96LNM} and \cite{S97AMS} for the precise statements and proofs in the Hamiltonian isotropic context, \cite{S01} for the volume preserving context with $p=2$, and \cite{S95} and \cite[\S~4.1.2]{BHS96LNM} for the reversible context~1 (see also a review in \cite{S98}). In other theorems, one employs the complicated technique of normal (Birkhoff-like) forms around an invariant torus (cf.\ \cite{CP11JMD,P13}), and the corresponding nondegeneracy conditions involve the coefficients of nonlinear terms of those normal forms, like the low order Taylor coefficients of the functions $F^\nu_j$ in~\eqref{Birk} at $0$. Such theorems are presented in e.g.\ \cite{JV97,JV01} (as well as in \cite{A62,A63,F04,CP09,CP10,CP11IM,F13,P15,Pnow}) for the Hamiltonian isotropic context and in \cite{S93} for the reversible context~1. The paper \cite{S97RCD} surveyed the theory of excitation of elliptic normal modes in Hamiltonian systems as it stood in 1997. This topic is also reviewed in \cite[p.~303]{AKN06}, \cite[Section~8.4]{BS10}, and \cite{S03}.

Within the Bryuno theory, one studies \emph{analytic} families of quasi-periodic invariant tori of dimensions $N\geq n$ passing through a given invariant $n$-torus, see e.g.\ \cite[Part~II]{B89}. Such families are subfamilies of Whitney smooth Cantor-like families of quasi-periodic invariant $N$-tori.

\section{Peculiarities of the Reversible Context~2}
\label{peculiarities}

Excitation of elliptic normal modes in the reversible $(n,a,b,s)$ context~2 with $b>a\geq 1$ (and $1\leq r\leq a$) has been entirely unexplored yet. In this context, the invariant $(n+r)$-tori $\fT$ can either be still within the scope of the reversible context~2 (if $b'=b-r>a$, i.e., if $r<b-a$) or pertain to the reversible context~1 (if $b\leq 2a$ and $b-a\leq r\leq a$ whence $b'=b-r\leq a$). However, there is some evidence that excitation of elliptic normal modes in the reversible context~2 is in fact impossible or at least drastically different from that in the three ``conventional'' KAM contexts discussed in Section~\ref{excitation}.

First, the number of zero Floquet exponents of the initial $n$-tori $\cT^\nu$ in the reversible context~2 is equal to $b-a$ (see Table~2), but the number $\kappa_1$ of zero Floquet exponents of the $(n+r)$-tori $\fT$ (if $r$ pairs of purely imaginary Floquet exponents of $\cT^\nu$ really excite) is \emph{smaller} than $\kappa_2=(b-a)+r$. Indeed, $\kappa_1=|b'-a|=|b-r-a|$. The absolute value of the difference of two positive numbers $b-a$ and $r$ is smaller than $\max(b-a,r)$, not to mention $b-a+r$. The ``defect'' $\kappa_2-\kappa_1$ is equal to
\[
\kappa_2-\kappa_1=b-a+r-|b-a-r|=2\min(b-a,r).
\]
If $r<2(b-a)$ then $\kappa_1<b-a$.

Second, let the $2a$ nonzero Floquet exponents of the initial $n$-tori $\cT^\nu$ be
\[
\pm i\beta^\nu_1,\ldots,\pm i\beta^\nu_r, \pm\lambda^\nu_1,\ldots,\pm\lambda^\nu_{a-r}.
\]
Then among the $2\min(a,b-r)$ nonzero Floquet exponents of the $(n+r)$-tori $\fT$, there are $2(a-r)$ Floquet exponents close to $\pm\lambda^\nu_1,\ldots,\pm\lambda^\nu_{a-r}$. But $a-r<\min(a,b-r)$. Consequently, among the $2\min(a,b-r)$ nonzero Floquet exponents of the $(n+r)$-tori $\fT$, there are also
\[
d=\min(a,b-r)-a+r=\min(r,b-a)=\frac{\kappa_2-\kappa_1}{2}\geq 1
\]
\emph{``new''} pairs, say, $\pm\chi_1,\ldots,\pm\chi_d$. Very roughly, ``symbolically'', and ``speculatively'' speaking, these $d$ pairs stem from ``coupling'' of $d$ zeroes associated with $d$ ``units'' coming from an $r$-parameter family of tori like~\eqref{C1Cr} and $d$ zeroes associated with $d$ ``minus units'' coming from
\[
b-a=(\codim\Fix G-\dim\cT^\nu)-\dim\Fix G.
\]
If some of these $2d$ ``new'' nonzero Floquet exponents are purely imaginary for some tori $\fT$, one may expect ``secondary'' excitation of elliptic normal modes of such tori $\fT$. In turn, ``secondary'' excitation may in principle be accompanied by ``tertiary'' one, and so on.

Third, the family of the invariant $(n+r)$-tori $\fT$ cannot be described ``pictorially'' by the assertion that each $n$-torus $\cT^\nu$ ``emits'' an $r$-parameter subfamily of tori $\fT$ invariant under the flow of a certain \emph{individual} vector field. Indeed, the whole family of tori $\fT$ is characterized by the difference $c'-s=a-(b-r)=r-(b-a)$, see Table~1. If $r\leq b-a$ then to each individual vector field, there corresponds no more than one of the tori $\fT$. But even if $r>b-a$, each individual vector field admits a Whitney smooth Cantor-like family of tori $\fT$ with the number of parameters equal to $r-(b-a)<r$.

Fourth, the tori $\fT$ invariant under the flows of $V_\mu$ are not necessarily concentrated around the family of the invariant $n$-tori $\cT^\nu$.

All these peculiarities of excitation of elliptic normal modes in the reversible context~2 can be easily seen in the following model (even ``toy'') example. Consider the systems
\begin{equation}
\dot{y}=u_\mu(z\barz), \quad \dot{z}=izW_\mu(z\barz)+zyv_\mu(z\barz)
\label{onetoy}
\end{equation}
dependent on an external parameter $\mu$ ranging in an open domain $P\subset\mR^s$ ($s\geq 1$), where $y$ and $z$ are the phase space variables ranging in a neighborhood $\cO_1(0)\subset\mR$ of the origin in $\mR$ and in a neighborhood $\cO_1^{\mC}(0)\subset\mC$ of the origin in $\mC$, respectively, while $u_\mu$, $v_\mu$, $W_\mu$ are smooth real-valued functions. We assume that $W_\mu(z\barz)>0$ for all $z\in\cO_1^{\mC}(0)$ and $\mu\in P$. The systems~\eqref{onetoy} are reversible with respect to the involution $G:(y,z)\mapsto(-y,\barz)$. The fixed point manifold of $G$ is the line $\Fix G=\{y=0, \: z\in\mR\}$.

If $u_\mu(0)=0$ then the system~\eqref{onetoy} at this value of $\mu$ admits exactly one equilibrium on $\Fix G$, namely, $(y=0, \: z=0)$. If $u_\mu(0)\neq 0$ then the system~\eqref{onetoy} has no such equilibria at all. Generically the condition $u_\mu(0)=0$ determines a hypersurface in $P$. The equilibrium $(y=0, \: z=0)$ pertains to the reversible $(0,1,2,s)$ context~2 with the Floquet exponents $0,\pm i\beta(\mu)$ where $\beta(\mu)=W_\mu(0)$. In the notation of Tables~1 and~2, this equilibrium is characterized by the quantities
\[
n=0, \;\; a=1, \;\; b=2, \;\; m=3, \;\; c=s-1, \;\; g=2.
\]

To watch excitation of the Floquet exponents $\pm i\beta(\mu)$, set $z=\sqrt{\rho}\,e^{i\varphi}$ and rewrite the systems~\eqref{onetoy} and the involution $G$ in the real variables $(y,\rho,\varphi)$:
\begin{gather}
\dot{y}=u_\mu(\rho), \quad \dot{\rho}=2\rho yv_\mu(\rho), \quad \dot{\varphi}=W_\mu(\rho);
\label{twotoy} \\
G:(y,\rho,\varphi)\mapsto(-y,\rho,-\varphi).
\nonumber
\end{gather}
The subsystems
\begin{equation}
\dot{y}=u_\mu(\rho), \quad \dot{\rho}=2\rho yv_\mu(\rho)
\label{toy}
\end{equation}
reversible with respect to the involution $\cG:(y,\rho)\mapsto(-y,\rho)$ do not depend on the phase $\varphi$ and possess the first integrals
\[
y^2 - \int^\rho\frac{u_\mu(\eta)}{\eta v_\mu(\eta)}\,d\eta.
\]
To each equilibrium $\fE(\mu,\rho_0)$ of the system~\eqref{toy} of the form $(y=0, \: \rho=\rho_0>0)\in\Fix\cG$, there corresponds a $G$-invariant cycle $\fT(\mu,\rho_0) = \{y=0, \: \rho=\rho_0\}$ of the system~\eqref{twotoy} with the frequency $W_\mu(\rho_0)$ close to $\beta(\mu)$ and with Floquet exponents which are generically nonzero (in the notation above, $r=d=1$). The cycles $\fT(\mu,\rho_0)$ pertain to the reversible $(1,1,1,s)$ context~1 and are characterized by the quantities
\[
n'=1, \;\; a=1, \;\; b'=1, \;\; m'=2, \;\; c'=s, \;\; g'=2.
\]

The equilibria $\fE(\mu,\rho_0)$ are determined by the equation $u_\mu(\rho_0)=0$. Generically for each $\mu$ in an open subdomain of $P$, there are several isolated equilibria of the form $(y=0, \: \rho=\rho_0>0)$, some of them being centers while the others being saddles. These equilibria do not accumulate to the point $(y=0, \: \rho=0)$, so that the cycles $\fT(\mu,\rho_0)$ do not accumulate to the point $(y=0, \: z=0)$ and are isolated in the phase space. The matrix of the linearization of the system~\eqref{toy} at an equilibrium $\fE(\mu,\rho_0)$ is
\begin{equation}
\begin{pmatrix} 0 & u'_\mu(\rho_0) \\
2\rho_0v_\mu(\rho_0) & 0 \end{pmatrix} \in \slsl(2,\mR),
\label{matrix}
\end{equation}
where $u'_\mu(\rho_0) = du_\mu(\rho)/d\rho|_{\rho=\rho_0}$. The variables $x=\varphi$, $X=(y,\rho-\rho_0)$ are Floquet coordinates for the cycle $\fT(\mu,\rho_0)$ of the system~\eqref{twotoy}, cf.\ the formulas~\eqref{Floq}, the matrix~\eqref{matrix} is the Floquet matrix of $\fT(\mu,\rho_0)$, and the eigenvalues
\[
\pm\chi(\mu,\rho_0) = \pm\bigl[ 2\rho_0u'_\mu(\rho_0)v_\mu(\rho_0) \bigr]^{1/2}
\]
of~\eqref{matrix} are the Floquet exponents of $\fT(\mu,\rho_0)$. If $u'_\mu(\rho_0)v_\mu(\rho_0)>0$ then $\fE(\mu,\rho_0)$ is a saddle of the system~\eqref{toy} and the Floquet exponents $\pm\chi(\mu,\rho_0)$ of $\fT(\mu,\rho_0)$ are real. If $u'_\mu(\rho_0)v_\mu(\rho_0)<0$ then $\fE(\mu,\rho_0)$ is a center and the Floquet exponents $\pm\chi(\mu,\rho_0)$ of $\fT(\mu,\rho_0)$ are purely imaginary.

Each center $\fE(\mu,\rho_0)$ is surrounded by $\cG$-invariant periodic trajectories of the system~\eqref{toy}. To these periodic trajectories, there correspond $G$-invariant $2$-tori of the system~\eqref{twotoy} which surround the cycle $\fT(\mu,\rho_0)$. Most of these tori carry quasi-periodic motions with strongly incommensurable frequencies. Such tori pertain to the reversible $(2,1,0,s)$ context~1, are characterized by the quantities
\[
n''=2, \;\; a=1, \;\; b''=0, \;\; m''=1, \;\; c''=s+1, \;\; g''=0,
\]
and can be regarded as a manifestation of ``secondary'' excitation of the Floquet exponents $\pm\chi(\mu,\rho_0)$ of $\fT(\mu,\rho_0)$.

\end{document}